\documentclass[11pt]{article}
\title{Stochastic PDEs with heavy-tailed noise}

\author{
Carsten Chong\thanks{Center for Mathematical Sciences, Technical University of Munich, Boltzmannstra\ss e 3, 85748 Garching, Germany, e-mail: carsten.chong@tum.de, url: www.statistics.ma.tum.de}
}

\usepackage{latexsym}
\usepackage{amssymb}
\usepackage{amsmath}
\usepackage{amsthm}
\usepackage{mathrsfs}
\usepackage[round]{natbib}
\usepackage{stmaryrd}
\usepackage{url}
\usepackage{longtable}
\usepackage{color}
\usepackage{dsfont}
\usepackage[T1]{fontenc}
\usepackage{lmodern}
\usepackage[linktocpage]{hyperref}

\usepackage{graphicx}
\usepackage[latin1]{inputenc}
\usepackage{enumitem}
\setenumerate{leftmargin=*}
\usepackage{amsopn}
\usepackage{comment}

\newcommand{\bfi}{\begin{fig}}
\newcommand{\efi}{\end{fig}}
\newcommand{\btab}{\begin{tab}}
\newcommand{\etab}{\end{tab}}
\newcommand{\barr}{\begin{array}}
\newcommand{\earr}{\end{array}}
\newcommand{\beq}{\begin{equation}}
\newcommand{\eeq}{\end{equation}}
\newcommand{\bdis}{\begin{displaymath}}
\newcommand{\edis}{\end{displaymath}\noindent}

\newcommand{\bbn}{\mathbb{N}}

\newcommand{\bbr}{\mathbb{R}}

\newcommand{\bbe}{\mathbb{E}}

\newcommand{\bbp}{\mathbb{P}}

\newcommand{\bbf}{\mathbb{F}}

\newcommand{\bone}{\mathds 1}

\newcommand{\halmos}{\quad\hfill $\Box$}

\newcommand{\cals}{{\cal S}}

\newcommand{\calf}{{\cal F}}

\newcommand{\cale}{{\cal E}}
\newcommand{\calp}{{\cal P}}
\newcommand{\calb}{{\cal B}}

\newcommand{\caln}{{\cal N}}

\newcommand{\calL}{{\cal L}}

\newcommand{\al}{{\alpha}}
\newcommand{\la}{{\lambda}}
\newcommand{\La}{{\Lambda}}

\newcommand{\ga}{{\gamma}}
\newcommand{\Ga}{{\Gamma}}

\newcommand{\si}{{\sigma}}
\newcommand{\om}{{\omega}}
\newcommand{\Om}{{\Omega}}

\newcommand{\ov}{\overline}
\newcommand{\un}{\underline}

\newcommand{\llb}{\llbracket}
\newcommand{\rrb}{\rrbracket}

\newcommand{\dd}{\mathrm{d}}
\newcommand{\ee}{\mathrm{e}}

\newcommand{\bb}{\mathrm{b}}

\newcommand{\loc}{\mathrm{loc}}

\newcommand{\Leb}{\mathrm{Leb}}

\newcommand{\pf}{\mathfrak{p}}
\newcommand{\qf}{\mathfrak{q}}

\makeatletter
\newcommand{\opnorm}{\@ifstar\@opnorms\@opnorm}
\newcommand{\@opnorms}[1]{%
  \left|\mkern-1.5mu\left|\mkern-1.5mu\left|
   #1
  \right|\mkern-1.5mu\right|\mkern-1.5mu\right|
}
\newcommand{\@opnorm}[2][]{%
  \mathopen{#1|\mkern-1.5mu#1|\mkern-1.5mu#1|}
  #2
  \mathclose{#1|\mkern-1.5mu#1|\mkern-1.5mu#1|}
}
\makeatother

\newtheoremstyle{neu}
    {11pt}      
    {11pt}      
    {}                  
    {}          
    {\bfseries} 
    {}          
    {1em}  
    {\textbf{\thmname{#1}\thmnumber{ #2}\thmnote{ (#3)}}}          
    
\newtheoremstyle{proof}
    {11pt}      
    {11pt}      
    {}                  
    {}          
    {\bfseries} 
    {}            
    {1em}          
    {\textbf{\thmname{#1}\thmnote{ #3}.}}          

\newtheorem{Theorem}{Theorem}[section]
\newtheorem{Corollary}[Theorem]{Corollary}
\newtheorem{Lemma}[Theorem]{Lemma}
\newtheorem{Proposition}[Theorem]{Proposition}
\theoremstyle{neu}
\newtheorem{Definition}[Theorem]{Definition}
\newtheorem{Example}[Theorem]{Example}
\newtheorem{Remark}[Theorem]{Remark}
\newtheorem{Assumption}{Assumption}

\theoremstyle{proof}
\newtheorem{Proof}{Proof}

\newcommand{\bthm}{\begin{Theorem}}
\newcommand{\ethm}{\end{Theorem}}

\newcommand{\bcor}{\begin{Corollary}}
\newcommand{\ecor}{\end{Corollary}}

\newcommand{\blem}{\begin{Lemma}}
\newcommand{\elem}{\end{Lemma}}

\newcommand{\bprop}{\begin{Proposition}}
\newcommand{\eprop}{\end{Proposition}}

\newcommand{\bdf}{\begin{Definition}}
\newcommand{\edf}{\end{Definition}}

\newcommand{\bex}{\begin{Example}}
\newcommand{\eex}{\end{Example}}

\newcommand{\brem}{\begin{Remark}}
\newcommand{\erem}{\end{Remark}}

\newcommand{\bass}{\begin{Assumption}}
\newcommand{\eass}{\end{Assumption}}

\newcommand{\bpr}{\begin{Proof}}
\newcommand{\epr}{\end{Proof}}

\newcommand{\benu}{\begin{enumerate}}
\newcommand{\eenu}{\end{enumerate}}

\newcommand{\bit}{\begin{itemize}}
\newcommand{\eit}{\end{itemize}}

\numberwithin{equation}{section}

\allowdisplaybreaks[2]

\begin{document}

\date{}

\maketitle

\begin{abstract}
	We analyze the nonlinear stochastic heat equation driven by heavy-tailed noise on unbounded domains and in arbitrary dimension. The existence of a solution is proved even if the noise only has moments up to an order strictly smaller than its Blumenthal-Getoor index. In particular, this includes all stable noises with index $\alpha<1+2/d$. Although we cannot show uniqueness, the constructed solution is natural in the sense that it is the limit of the solutions to approximative equations obtained by truncating the big jumps of the noise or by restricting its support to a compact set in space. Under growth conditions on the nonlinear term we can further derive moment estimates of the solution, uniformly in space. Finally, the techniques are shown to apply to Volterra equations with kernels bounded by generalized Gaussian densities. This includes, for instance, a large class of uniformly parabolic stochastic PDEs.
\end{abstract}

\vfill

\noindent
\begin{tabbing}
\> secondary: \,\,\,60G52, 60G48, 60G51, 60G57
\end{tabbing}

\vspace{1cm}


\vspace{0.5cm}

\newpage

\section{Introduction}\label{intro}

In this paper we investigate the nonlinear \emph{stochastic heat equation}
\begin{align} \partial_t Y(t,x) &= \Delta Y(t,x) + \si(Y(t,x))\dot M(t,x),\quad (t,x)\in\bbr_+\times\bbr^d,\nonumber\\
 Y(0,x)&=\psi(x),	\label{SPDE}
\end{align}
where $M$ is some random noise, the function $\si$ is globally Lipschitz and $\psi$ is some given initial condition. In the theory of stochastic partial differential equations (stochastic PDEs) there are various ways to make sense of the formal equation \eqref{SPDE}. We refer to \cite{Peszat07}, \cite{Liu15} and \cite{Holden10} for detailed accounts within the semigroup, the variational and the white noise approach, respectively. Our analysis is based on the random field approach going back to \cite{Walsh86} where the differential equation \eqref{SPDE} is interpreted as the integral equation
\beq\label{SPDE-heat} Y(t,x)=Y_0(t,x)+\int_0^t \int_{\bbr^d} g(t-s,x-y)\si(Y(s,y))\,M(\dd s,\dd y),\quad (t,x)\in \bbr_+\times\bbr^d,\eeq
where $g$ is the heat kernel
\beq\label{hk} g(t,x)=\frac{1}{(4\pi t)^{d/2}} \exp\left(-\frac{|x|^2}{4t}\right) \bone_{(0,\infty)}(t),\quad (t,x)\in\bbr_+\times\bbr^d,\eeq
$Y_0$ is related to the initial condition $\psi$ via
\[ Y_0(t,x) = \int_{\bbr^d} g(t,x-y)\psi(y)\,\dd y,\quad (t,x)\in\bbr_+\times\bbr^d, \]
and the integral with respect to $M$ is understood in It\^o's sense, see e.g. \cite{Bichteler83}, \cite{Walsh86} and \cite{Chong15}. 

There are multiple reasons for the broad interest in studying the stochastic heat equation \eqref{SPDE}: on the one hand, it is of great significance within mathematics because, for example, it has strong connections to branching interacting particles (see e.g. \cite{Carmona94}, \cite{Dawson93} and \cite{Mueller15}) and it arises from the famous Kardar-Parisi-Zhang equation via the Cole-Hopf transformation \citep{Corwin12, Hairer13}. On the other hand, the stochastic heat equation has found applications in various disciplines like turbulence modeling \citep{Davies04, BN11-2}, astrophysics \citep{Jones99}, phytoplankton modeling \citep{ElSaadi15} and neurophysiology \citep{Walsh81, Tuckwell83}.

If $M$ is a Gaussian noise, the initial theory developed by \cite{Walsh86} has been comprehensively extended and many properties of the solution to the stochastic heat equation or variants hereof have been established, see \cite{Khoshnevisan14}, for example. When it comes to non-Gaussian noise, the available literature on stochastic PDEs within the random field approach is much less. While the papers by \cite{Albeverio98} and \cite{Applebaum00} remain in the $L^2$-framework of \cite{Walsh86}, \cite{SLB98} is one of the first to treat L\'evy-driven stochastic PDEs in $L^p$-spaces with $p<2$. The results are extended in \cite{Chong16} to Volterra equations with L\'evy noise, on finite as well as on infinite time domains. In both \cite{SLB98} and \cite{Chong16}, one crucial assumption the L\'evy noise has to meet is that its L\'evy measure, say $\la$, must satisfy
\beq\label{intlevy} \int_\bbr |z|^p \,\la(\dd z)<\infty \eeq
for some $p<1+2/d$. This a priori excludes stable noises, or more generally, any noise with less moments than its Blumenthal-Getoor index. Of course, if the noise intensity decreases sufficiently fast in space, for example, if it only acts on a bounded domain instead of the whole $\bbr^d$, this problem can be solved by localization, see \cite{Balan14} and \cite{Chong16}. Also on a bounded space domain, \cite{Mueller98} considers the stochastic heat equation with spectrally positive stable noise of index smaller than one and the non-Lipschitz function $\si(x)=x^\ga$ ($\ga<1$). To our best knowledge, the paper by \cite{Mytnik02}, where the same function $\si$ is considered with spectrally positive stable noises of index bigger than one, is the only work that investigates the stochastic heat equation on the whole $\bbr^d$ with noises violating \eqref{intlevy} for all $p\geq 0$. The purpose of this article is to present results for Lipschitz $\si$ and more general heavy-tailed noises $M$, including for all dimensions $d$ stable noises of index $\al<1+2/d$ (with no drift in the case $\al<1$). 

More specifically, if $M$ is a homogeneous pure-jump L\'evy basis with L\'evy measure $\la$ and $\psi$ is, say, a bounded deterministic function, we will prove in Theorem~\ref{main} the existence of a solution to \eqref{SPDE-heat} under the assumption that
\[ \int_{|z|\leq 1} |z|^p \,\la(\dd z) + \int_{|z|>1} |z|^q \,\la(\dd z) < \infty \]
for some $p<1+2/d$ and $q>p/(1+(1+2/d-p))$. Although the proof does not yield uniqueness in a suitable $L^p$-space, we can show in Theorem~\ref{approx} that our solution is the limit in probability of the solutions to \eqref{SPDE-heat} when either the jumps of $M$ are truncated at increasing levels or when $M$ is restricted to increasing compact subsets of $\bbr^d$. Furthermore, if $|\si(x)|\leq C(1+|x|^\ga)$ with an exponent $\ga<1\wedge q/p$, the constructed solution has $q$-th order moments that are uniformly bounded in space, see Theorem~\ref{moment}. Finally, in Theorem~\ref{extension}, we extend the previous results to Volterra equations with kernels bounded by generalized Gaussian densities, which includes uniformly parabolic stochastic PDEs as particular examples.

\section{Preliminaries}\label{prelim}
Underlying the whole paper is a stochastic basis $(\Om,\calf,\bbf=(\calf_t)_{t\in \bbr_+},\bbp)$ satisfying the usual hypotheses of right-continuity and completeness that is large enough to support all considered random elements. For $d\in\bbn$ we equip $\tilde\Om := \Om\times \bbr_+ \times\bbr^d$ with the \emph{tempo-spatial predictable $\si$-field} $\tilde\calp:=\calp\otimes\calb(\bbr^d)$ where $\calp$ is the usual predictable $\si$-field and $\calb(\bbr^d)$ is the Borel $\si$-field on $\bbr^d$. With a slight abuse of notation, $\tilde\calp$ also denotes the space of \emph{predictable}, that is, $\tilde\calp$-measurable processes $\tilde\Om\to\bbr$. Moreover, $\tilde\calp_\bb$ is the collection of all $A\in\tilde\calp$ satisfying $A\subseteq \Om\times[0,k]\times[-k,k]^d$ for some $k\in\bbn$. Next, writing $p^\ast:=p\vee 1$ and $p_\ast:=p\wedge1$ for $p\in[0,\infty)$, we equip the space $L^p=L^p(\Om,\calf,\bbp)$ with the usual topology induced by $\|X\|_{L^p} := \bbe[|X|^p]^{1/p^\ast}$ for $p>0$ and $\|X\|_{L^0} := \bbe[|X|\wedge1]$ for $p=0$. Further notations and abbreviations include $\llbracket R,S \rrbracket:=\{(\om,t)\in\Om\times \bbr_+\colon R(\om)\leq t\leq S(\om)\}$ for two $\bbf$-stopping times $R$ and $S$ (other stochastic intervals are defined analogously), $|\mu|$ for the total variation measure of a signed Borel measure $\mu$, $x^{(i)}$ for the $i$-th coordinate of a point $x\in\bbr^d$ and $x_+:=x\vee 0$ for $x\in\bbr$. Last but not least, the value of all constants $C$ and $C_T$ may vary from line to line without further notice.

In this paper the noise $M$ that drives the equations \eqref{SPDE} and \eqref{SPDE-heat} will be assumed to have the form
\begin{align} M(\dd t,\dd x) &= b(t,x)\,\dd(t,x) + \rho(t,x)\,W(\dd t,\dd x) + \int_E \un\delta(t,x,z)\,(\pf-\qf)(\dd t,\dd x,\dd z) \nonumber\\
&\quad + \int_E \ov\delta(t,x,z)\,\pf(\dd t,\dd x,\dd z) \label{candecM} \end{align}
with the following specifications:
\bit
 \item $(E,\cale)$ is an arbitrary Polish space equipped with its Borel $\si$-field,
 \item $b,\rho\in\tilde\calp$,
 \item $\delta:=\un\delta+\ov\delta:=\delta\bone_{\{|\delta|\leq 1\}} + \delta\bone_{\{|\delta|>1\}}$ is a $\tilde\calp\otimes\cale$-measurable function,
 \item $W$ is a Gaussian space--time white noise relative to $\bbf$ (see \cite{Walsh86}),
 \item $\pf$ is a homogeneous Poisson random measure on $\bbr_+\times\bbr^d\times E$ relative to $\bbf$ (see Definition~II.1.20 in \cite{Jacod03}), whose intensity measure disintegrates as $\qf(\dd t,\dd x,\dd z) = \dd t\,\dd x\,\la(\dd z)$ with some $\si$-finite measure $\la$ on $(E,\cale)$.
\eit
Moreover, all coefficients are assumed to be ``locally integrable'' in the sense that the random variable 
\begin{align*} M(A)&=\int_{\bbr_+\times\bbr^d} \bone_A(t,x)b(t,x)\,\dd(t,x) + \int_{\bbr_+\times\bbr^d} \bone_A(t,x)\rho(t,x)\,W(\dd t,\dd x)\\
&\quad + \int_{\bbr_+\times\bbr^d\times E} \bone_A(t,x)\un\delta(t,x,z)\,(\pf-\qf)(\dd t,\dd x,\dd z)\\
&\quad + \int_{\bbr_+\times\bbr^d\times E} \bone_A(t,x)\ov\delta(t,x,z)\,\pf(\dd t,\dd x,\dd z)  
\end{align*}
is well defined for all $A\in\tilde\calp_\bb$. In analogy to the notion of It\^o semimartingales in the purely temporal case, we call the measure $M$ in \eqref{candecM} an \emph{It\^o basis}. We make the following structural assumptions on $M$: there exist exponents $p, q \in (0,2]$ (without loss of generality we assume that $q\leq p$), numbers $\beta_N \in \bbr_+$, $\bbf$-stopping times $\tau_N$, and deterministic positive measurable functions $j_N(z)$ such that
	\benu
		\item[M1.] $\tau_N>0$ for all $N\in\bbn$ and $\tau_N\uparrow \infty$ a.s.,
		\item[M2.] $\int_E j_N(z)\,\la(\dd z)<\infty$,
		\item[M3.] $|b(\om,t,x)|, |\rho^2(\om,t,x)|\leq \beta_N$ and $|\un\delta(\om,t,x,z)|^p + |\ov\delta(\om,t,x,z)|^q \leq j_N(z)$ for all $(\om,t,x)\in\tilde\Om$ with $t\leq \tau_N(\om)$ and $z\in E$.
	\eenu
If $p\leq 1$ (resp. $q\geq1$), we can a.s. define $b_0(t,x):=b(t,x)-\int_E \un\delta(t,x,z) \,\la(\dd z)$ (resp. $b_1(t,x):=b(t,x) + \int_E \ov\delta(t,x,z)\,\la(\dd z)$) for $(t,x)\in\bbr_+\times\bbr^d$. In this case, we assume without loss of generality that $|b_0(\om,t,x)|$ (resp. $|b_1(\om,t,x)|$) is bounded by $\beta_N$ as well when $t\leq \tau_N(\om)$.

\bex\label{Levybasis} If $M$ is an It\^o basis with $E=\bbr$, deterministic coefficients $b(t,x)=b$, $\rho(t,x)=\rho$ and $\delta(t,x,z)=z$, and $\la$ satisfies $\la(\{0\})=0$ and $\int_\bbr (1\wedge z^2)\,\la(\dd z)<\infty$, then $M$ is a \emph{homogeneous L\'evy basis} (and $\dot M$ a space--time white noise) with L\'evy measure $\la$. In this case, conditions M1--M3 reduce to the requirement that
\beq\label{Levypq} \int_{|z|\leq 1} |z|^p\,\la(\dd z) + \int_{|z|>1} |z|^q\,\la(\dd z)<\infty,\eeq
with $\tau_N:=+\infty$ and $j_N(z):=|z|^p\bone_{[-1,1]}(z)+|z|^q\bone_{\bbr\setminus[-1,1]}(z)$. If $p\leq 1$ (resp. $q\geq1$), $b_0$ (resp. $b_1$) is the drift (resp. the mean) of $M$.\halmos 
\eex


Regarding integration with respect to It\^o bases, essentially, as all integrals are taken on finite time intervals, the reader only has to be acquainted with the integration theory with respect to Gaussian white noise (see \cite{Walsh86}) and (compensated) Poisson random measures (see Chapter~II of \cite{Jacod03}). But at one point in the proof of Theorem~\ref{main} below, we need the more abstract theory behind. So we briefly recall this, all omitted details can be found in e.g. \cite{Bichteler83} or \cite{Chong15}. The stochastic integral of a simple function is defined in the canonical way: if $S(t,x)=\sum_{i=1}^r a_i\bone_{A_i}$ with $r\in\bbn$, $a_i\in\bbr$ and $A_i\in\tilde\calp_\bb$, then 
\[ \int S\,\dd M := \int S(t,x) \,M(\dd t,\dd x) := \sum_{i=1}^r a_i M(A_i). \]
In order to extend this integral from the class of simple functions $\cals$ to general predictable processes, it is customary to introduce the \emph{Daniell mean}
\[ \|H\|_{M,p} := \sup_{S\in\cals,|S|\leq |H|} \left\|\int S\, \dd M\right\|_{L^p}\]
for $H\in\tilde\calp$ and $p\in\bbr_+$. A predictable function $H$ is called \emph{$L^p$-integrable} (or simply \emph{integrable} if $p=0$) if there exists a sequence $(S_n)_{n\in\bbn}$ of simple functions in $\cals$ with $\|S_n-H\|_{M,p}\to0$ as $n\to\infty$. In this case, the \emph{stochastic integral}
\[ \int H \,\dd M := \int H(t,x)\,M(\dd t,\dd x) := \lim_{n\to\infty} \int S_n(t,x)\, M(\dd t,\dd x) \]
exists as a limit in $L^p$ and does not depend on the choice of $S_n$. It is convenient to indicate the domain of integration by $\int_A$ or $\int_0^t\int_{\bbr^d} = \int_{(0,t)\times\bbr^d}$, for example. In the latter case, contrary to usual practice, the endpoint $t$ is excluded. 

Although the Daniell mean seems to be an awkward expression, it can be computed or estimated effectively in two important situations: if $M$ is a strict random measure, that is, if it can be identified with a measure-valued random variable $M(\om,\cdot)$, then we simply have $\|H\|_{M,p}=\left\|\int |H|\,\dd|M|\right\|_{L^p}$; if $M$ is a local martingale measure, that is, for every $A\in\tilde\calp_\bb$ the process $t\mapsto M(A\cap(\Om\times[0,t]\times\bbr^d))$ has a version that is an $\bbf$-local martingale, then by the Burkholder-Davis-Gundy inequalities (in conjunction with Theorem~VII.104 of \cite{Dellacherie82}) there exists, for every $p\in[1,\infty)$, a constant $C_p\in\bbr_+$ such that
\[ \|H\|_{M,p} \leq C_p \left\|\left(\int H^2\,\dd[M]\right)^{1/2} \right\|_{L^p}  \]
for all $H\in\tilde\calp$. Here $[M](\dd t,\dd x):=\rho^2(t,x)\,\dd(t,x) + \int_E \delta^2(t,x,z)\,\pf(\dd t,\dd x,\dd z)$ is the quadratic variation measure of $M$.

Finally, for $p\in(0,\infty)$, we define the space $B^p$ (resp. $B^p_\loc$) as the collection of all $\phi\in\tilde\calp$ for which we have
\begin{align*} \|\phi\|_{p,T} &:= \sup_{(t,x)\in [0,T] \times \bbr^d} \|\phi(t,x)\|_{L^p} <\infty\\
\Bigg(\text{resp.}\quad \|\phi\|_{p,T,R} &:= \sup_{(t,x)\in [0,T] \times [-R,R]^d} \|\phi(t,x)\|_{L^p} <\infty \Bigg)
\end{align*}
for all $T\in\bbr_+$ (resp. for all $T, R \in\bbr_+$). Moreover, if $\tau$ is an $\bbf$-stopping time, we write $\phi\in B^p(\tau)$ (resp. $\phi\in B^p_\loc(\tau)$) if $\phi\bone_{\llb 0, \tau\rrb}$ belongs to $B^p$ (resp. $B^p_\loc$). 

\section{Main results}\label{mainres}

In \cite{SLB98} and \cite{Chong16}, Equation~\eqref{SPDE-heat} driven by a L\'evy basis is proved to possess a solution in some space $B^p$ under quite general assumptions. The most restrictive condition, however, is that the L\'evy measure $\la$ has to satisfy 
\beq\label{pintegrab} \int_\bbr |z|^p\,\la(\dd z)< \infty.\eeq
Although there is freedom in choosing the value of $p\in(0,1+2/d]$, there is no possibility to choose this exponent for the small jumps (i.e., $|z|\leq 1$) and the big jumps (i.e., $|z|>1$) separately. So for instance, stable L\'evy bases are always excluded. The main goal of this paper is to describe a way of how one can, to a certain extent, choose different exponents $p$ in \eqref{pintegrab} for the small and the large jumps, respectively. In fact, we can consider a slightly more general equation than \eqref{SPDE-heat}, namely the \emph{stochastic Volterra equation}
\beq\label{SPDE-var} Y(t,x)=Y_0(t,x)+\int_0^t \int_{\bbr^d} G(t,x;s,y)\si(Y(s,y))\,M(\dd s,\dd y),\quad (t,x)\in \bbr_+\times\bbr^d,\eeq
under the following list of assumptions.
\benu
\item[A1.] For every $N\in\bbn$ we have that $Y_0\in B^p(\tau_N)$.
\item[A2.] There exists $C\in\bbr_+$ such that $|\si(x)-\si(y)|\leq C|x-y|$ for all $x, y \in \bbr$.
\item[A3.] The function $G\colon (\bbr_+\times\bbr^d)^2 \to \bbr$ is measurable and for all $T\in\bbr_+$ there exists $C_T\in\bbr_+$ with $|G(t,x;s,y)| \leq C_T g(t-s,x-y)$ for all $(t,x),(s,y)\in[0,T]\times\bbr^d$,
where $g$ is the heat kernel \eqref{hk}.
\item[A4.] The exponents $p$ and $q$ in M3 satisfy $0<p<1+2/d$ and $p/(1+(1+2/d-p))<q\leq p$.
\item[A5.] If $p<2$, we assume $\rho\equiv0$; if $p<1$, we also require that $b_0\equiv0$.
\eenu

As it is usual in the context of stochastic PDEs, we call $Y\in\tilde\calp$ a \emph{solution} to \eqref{SPDE-var} if for all $(t,x)\in \bbr_+\times\bbr^d$ the stochastic integral on the right-hand side of \eqref{SPDE-var} is well defined and the equation itself holds a.s. Two solutions $Y_1$ and $Y_2$ are identified if they are versions of each other.

\bthm\label{main} 
Let $M$ be an It\^o basis satisfying M1--M3. Then under A1--A5 there exists a sequence $(\tau(N))_{N\in\bbn}$ of $\bbf$-stopping times increasing to infinity a.s. such that Equation~\eqref{SPDE-var} has a solution $Y$ belonging to $B^p_\loc(\tau(N))$ for every $N\in\bbn$. One possible choice of $(\tau(N))_{N\in\bbn}$ is given in \eqref{tauN} below.
\ethm

Figure~\ref{fig} illustrates the possible choices for the exponents $p$ and $q$ such that Theorem~\ref{main} is applicable. For each dimension $d\in\bbn$, the exponent $p$ must be smaller than $(1+2/d)\wedge 2$ and the exponent $q$ larger than $p/(1+(1+2/d-p))$, which corresponds to the area above the bold lines in Figure~\ref{fig}. If $q\geq p$, the result is well known from \cite{SLB98} and \cite{Chong16}. So the new contribution of Theorem~\ref{main} pertains to the area above the bold lines and below the diagonal $q=p$. As we can see, all stable noises with index $\alpha<1+2/d$ (and zero drift if $\al<1$) are covered: they are ``located'' infinitesimally below the line $q=p$ in the figure. Moreover, we observe that there exists a non-empty region below the diagonal that constitutes valid parameter choices for all dimensions $d\in\bbn$: namely when $p\leq 1$ and $p/(2-p)<q\leq p$.

\begin{figure}
	\centering
	\includegraphics[height=8cm]{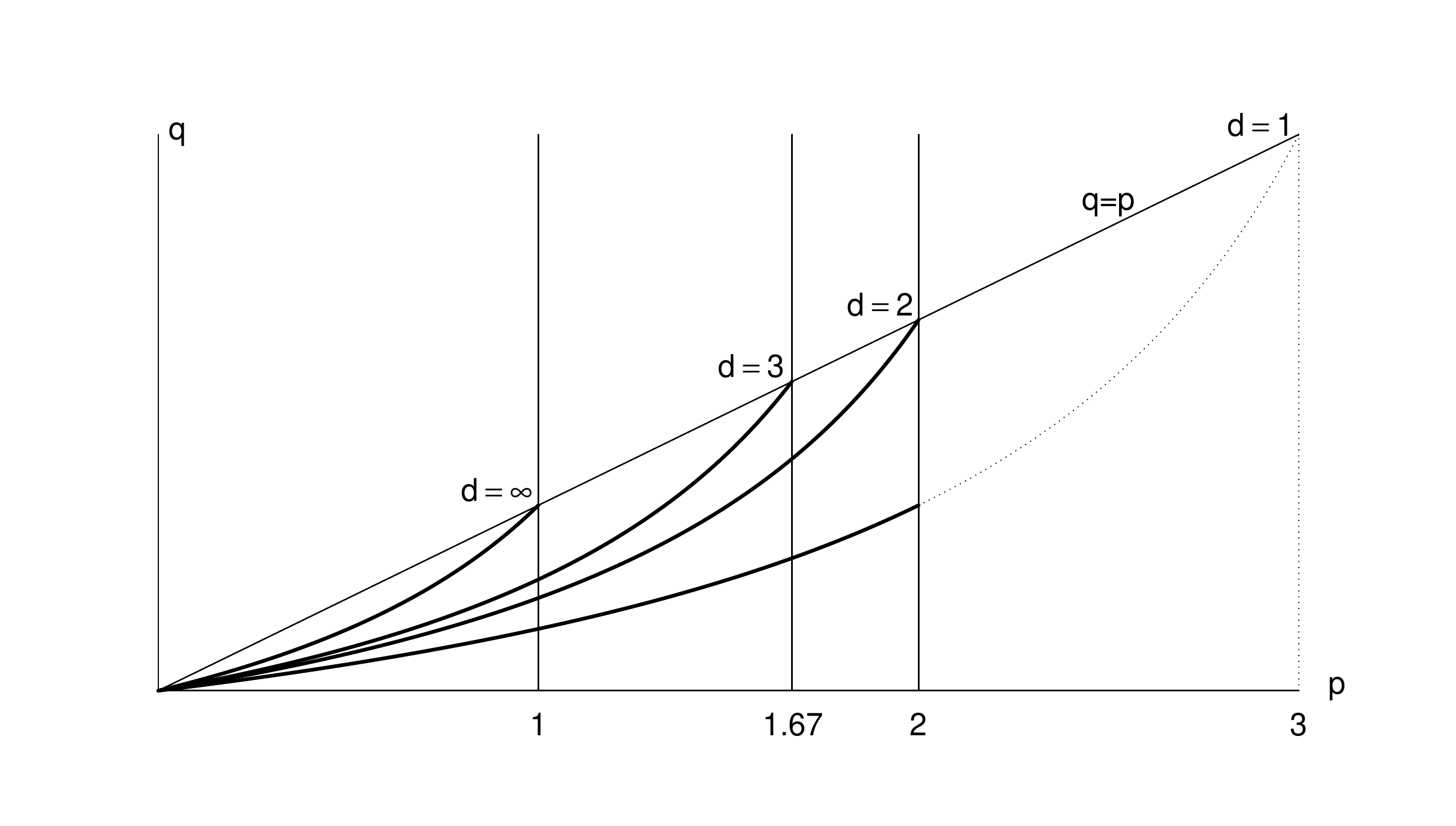}
	\vspace{-\baselineskip}
	\caption{Constraints for $p$ and $q$ in Theorem~\ref{main} dependent on the dimension $d$}\label{fig}
\end{figure}

Put in a nutshell, the method of \cite{SLB98} and \cite{Chong16} to construct a solution to \eqref{SPDE-var} is this: first, one defines the integral operator $J$ that maps $\phi\in\tilde\calp$ to
\beq\label{Jdef}J(\phi)(t,x):=Y_0(t,x)+\int_0^t\int_{\bbr^d} G(t,x;s,y)\si(\phi(s,y))\,M(\dd s,\dd y),\quad(t,x)\in \bbr_+\times\bbr^d, \eeq
when the right-hand side is well defined, and to $+\infty$ otherwise. Then one establishes moment estimates in such a way that $J$ becomes, at least locally in time, a contractive self-map on some $B^p$-space. Then the Banach fixed point theorem yields existence and uniqueness of solutions. 
However, if $q<p$, the problematic point is that the $p$-th order moment of $J(\phi)$ is typically infinite, while taking the $q$-th order moment, which would be finite, produces a non-contractive iteration. The idea in this article to circumvent this is therefore to consider a modified Picard iteration: we still take $p$-th moment estimates (in order to keep contractivity) but stop the processes under consideration before they become too large (in order to keep the $p$-th moments finite). The following lemma makes our stopping strategy precise.

\blem\label{tauNK} Let $h(x):=1+|x|^\eta$ with some $\eta\in\bbr_+$ and define for $N\in\bbn$ the stopping times 
\beq\label{tauN} \tau(N):=\inf\left\{T\in\bbr_+\colon \int_0^T\int_{\bbr^d} \int_E \bone_{\{|\ov\delta(t,x,z)|>N h(x)\}}\,\pf(\dd t,\dd x,\dd z) \neq 0\right\}\wedge \tau_N, \eeq
where $(\tau_N)_{N\in\bbn}$ are the stopping times in M1. 
If $q$ is the number in M3 and $\eta>d/q$, then a.s. $\tau(N)>0$ for all $N\in\bbn$ and $\tau(N)\uparrow\infty$ a.s. for $N\uparrow \infty$.
\elem

\bpr Define $U_0:=\emptyset$ and
\[ U_n:=\big\{x\in\bbr^d\colon |x|\leq \pi^{-1/2}\Gamma(1+d/2)^{1/d}n^{1/d}\big\}\setminus U_{n-1},\quad n\in\bbn. \]
In other words, $(U_n)_{n\in\bbn}$ forms a partition of $\bbr^d$ into concentric spherical shells of Lebesgue measure $1$. Then there exists a constant $C\in\bbr_+$ which is independent of $x$ and $n$ such that $h(x)\geq C (n-1)^{\eta/d}$ for all $x\in U_n$ and $n\in\bbn$. Next, observe that for all $T\in\bbr_+$ we a.s. have  
\begin{align*} &\int_0^T\int_{\bbr^d} \int_E \bone_{\{|\ov\delta(t,x,z)|>N h(x)\}}\bone_{\llb 0,\tau_N\rrb}(t)\,\pf(\dd t,\dd x,\dd z) \\
&\qquad\leq \int_0^T\int_{\bbr^d} \int_E \bone_{\{j_N(z)^{1/q}>Nh(x)\}}\,\pf(\dd t,\dd x,\dd z)\\
&\qquad\leq \sum_{n=1}^\infty \int_0^T\int_{\bbr^d} \int_E \bone_{U_n}(x) \bone_{\{j_N(z)^{1/q}>NC (n-1)^{\eta/d}\}}\,\pf(\dd t,\dd x,\dd z). \label{number}\end{align*}
If we show that the last expression is finite a.s., we can conclude, on the one hand, that $\tau(N)>0$ a.s. for all $N\in\bbn$. On the other hand, it also implies that $\tau(N) \uparrow \infty$ for $N\uparrow \infty$ because then, for any $T\in\bbr_+$,  the left-hand side of the previous display would be $0$ as soon as $N$ is large enough. So with $a_n:=N C (n-1)^{\eta/d}$ and using the Borel-Cantelli lemma, it remains to verify that
\beq\label{bc} \sum_{n=1}^\infty \bbp[\pf([0,T]\times U_n \times \{j_N(z)>a_n^q\})\neq0] < \infty. \eeq
Indeed,
\begin{align*} \bbp[\pf([0,T]\times U_n \times \{j_N(z)>a_n^q\})\neq0] &= 1-\ee^{-\qf([0,T]\times U_n \times \{j_N(z)>a_n^q\})}\\ 
&\leq \qf([0,T]\times U_n \times \{j_N(z)>a_n^q\})\\
&=T\la(\{j_N(z)>a_n^q\})\leq T \int_E j_N(z)\,\la(\dd z) a_n^{-q} .\end{align*}
The last term is of order $n^{-\eta q/d}$. Since $\eta>d/q$, it is summable and leads to \eqref{bc}. \halmos
\epr

The stopping times in the previous lemma are reminiscent of a well established technique of solving stochastic differential equations driven by semimartingales \citep{Protter05}. Here, if the jumps of the driving semimartingale lack good integrability properties, one can still solve the equation until the first big jump occurs, then continue until the next big jump occurs and so on, thereby obtaining a solution for all times. 
In principle, this is exactly the idea we employ in the proof of Theorem~\ref{main}. However, in our tempo-spatial setting, the law of $M$ is typically not equivalent to the law of its truncation at a fixed level on any set $[0,T]\times\bbr^d$: Infinitely many big jumps usually occur already immediately after time zero. This is the reason why in Lemma~\ref{tauNK} above we have to take a truncation level $h$ that increases sufficiently fast in the spatial coordinate.


As a next step we introduce for each $N\in\bbn$ a truncation of $M$ by
\begin{align} M^N(\dd t,\dd x) :=&~ \bone_{\llb 0, \tau_N\rrb} (t)\left(b(t,x)\,\dd(t,x) + \rho(t,x)\,W(\dd t,\dd x) + \int_E \un\delta(t,x,z)\,(\pf-\qf)(\dd t,\dd x,\dd z)\right.\nonumber\\
& +\left. \int_E \ov\delta(t,x,z)\bone_{\{|\delta(t,x,z)|\leq Nh(x)\}}\,\pf(\dd t,\dd x,\dd z)\right). \label{MNK}\end{align}
An important step in the proof of Theorem~\ref{main} is to prove that \eqref{SPDE-var} has a solution when $M$ is substituted by $M^{N}$. To this end, we establish several moment estimates.
\blem\label{mom} Let $N\in\bbn$ be arbitrary but fixed and $J_N$ the integral operator defined in the same way as $J$ in \eqref{Jdef} but with $M^{N}$ instead of $M$. 
Under M1--M3 and A1--A5 the following estimates hold.
\benu
\item For all $T\in\bbr_+$ there exists a constant $C_T\in\bbr_+$ such that for all $\phi\in\tilde\calp$ and $(t,x)\in[0,T]\times\bbr^d$ we have
\begin{align} &\|J_N(\phi)(t,x)\|_{L^p} \leq \|Y_0(t,x)\bone_{\llb 0, \tau_N\rrb}(t)\|_{L^p}\nonumber\\
&\quad\quad+ C_T \left(\int_0^t\int_{\bbr^d} g^p(t-s,x-y)(1+\|\phi(s,y)\|^{p^\ast}_{L^p})h(y)^{p-q}\,\dd(s,y)\right)^{\frac{1}{p^\ast}}\nonumber\\
&\quad\quad+ C_T \left(\int_0^t\int_{\bbr^d} g(t-s,x-y)(1+\|\phi(s,y)\|^p_{L^p})h(y)^{p-q}\,\dd(s,y)\right)^{\frac{1}{p}}\bone_{\{p\geq1\}}. \label{JN1}\end{align}
\item For all $T\in\bbr_+$ there exists a constant $C_T\in\bbr_+$ such that for all $(t,x)\in[0,T]\times\bbr^d$ and $\phi_1,\phi_2\in\tilde\calp$ with $J_N(\phi_1)(t,x), J_N(\phi_2)(t,x)<\infty$ we have
\begin{align*} &\|J_N(\phi_1)(t,x)-J_N(\phi_2)(t,x))\|_{L^p}  \\
&\quad\quad\leq C_T \left(\int_0^t\int_{\bbr^d} g^p(t-s,x-y)\|\phi_1(s,y)-\phi_2(s,y)\|^{p^\ast}_{L^p}h(y)^{p-q}\,\dd(s,y)\right)^{\frac{1}{p^\ast}} \\
&\quad\quad\quad+ C_T \left(\int_0^t\int_{\bbr^d} g(t-s,x-y)\|\phi_1(s,y)-\phi_2(s,y)\|^p_{L^p}h(y)^{p-q}\,\dd(s,y)\right)^{\frac{1}{p}}\bone_{\{p\geq1\}}. \end{align*}
\eenu
\elem

\bpr Both parts can be treated in a similar fashion. We only prove (1). Starting with the case $p\geq 1$, we use the fact that $\ov\delta = \delta\bone_{ \{| \delta|>1 \} }$, $\un\delta = \delta\bone_{ \{| \delta|\leq 1 \} }$ and $h \geq 1$ to decompose $M^{N}(\dd t,\dd x) = L^{N}(\dd t,\dd x) + B^{N}_1(\dd t,\dd x)$ where
\begin{align}
L^{N}(\dd t,\dd x)&:= \bone_{\llb 0, \tau_N\rrb} (t)\left(\rho(t,x)\,W(\dd t,\dd x) + \int_E \delta(t,x,z)\bone_{\{|\delta(t,x,z)|\leq Nh(x)\}}\,(\pf-\qf)(\dd t,\dd x,\dd z)\right), \nonumber\\
B^{N}_1(\dd t,\dd x)&:= \bone_{\llb 0, \tau_N\rrb} (t) \left(b(t,x)\,\dd(t,x)+\int_E \ov\delta(t,x,z)\bone_{\{|\delta(t,x,z)|\leq Nh(x)\}}\,\qf(\dd t,\dd x, \dd z)\right). \label{LBNK}
\end{align}
In the following $C_T$ denotes a positive constant independent of $\phi$ (but possibly dependent on $p$, $q$, $d$ and $N$). Recalling that $\rho=0$ unless $p=2$, using the Burkholder-Davis-Gundy inequalities and applying the inequality $(x+y)^r\leq x^r+y^r$ for $x,y\in\bbr_+$ and $r\in[0,1]$ to the Poisson integral, we deduce for all $(t,x)\in[0,T]\times\bbr^d$ that
\begin{align*}
&~\left\|\int_0^t\int_{\bbr^d} G(t,x;s,y)\si(\phi(s,y))\,L^{N}(\dd s,\dd y)\right\|_{L^p}\\
\leq&~ C_T \bbe\Bigg[\int_0^t\int_{\bbr^d} |g(t-s,x-y)\si(\phi(s,y))\rho(s,y)|^2\bone_{\llb 0, \tau_N\rrb}(s)\,\dd(s,y)\\
& + \left(\int_0^t\int_{\bbr^d}\int_E |g(t-s,x-y)\si(\phi(s,y))\delta(s,y,z)|^2\bone_{\{|\delta(s,y,z)|\leq Nh(y)\}}\bone_{\llb 0, \tau_N\rrb}(s)\,\pf(\dd s,\dd y,\dd z)\right)^{\frac{p}{2}}\Bigg]^{\frac{1}{p}}\\
\leq&~ C_T \bbe\Bigg[\int_0^t\int_{\bbr^d} g^2(t-s,x-y) (1+|\phi(s,y)|^2)\bone_{\{p=2\}}\,\dd(s,y)\\
& + \int_0^t\int_{\bbr^d}\int_E |g(t-s,x-y)\delta(s,y,z)|^p(1+|\phi(s,y)|^p) \bone_{\{|\delta(s,y,z)|\leq Nh(y)\}}\bone_{\llb 0, \tau_N\rrb}(s)\,\pf(\dd s,\dd y,\dd z)\Bigg]^{\frac{1}{p}}\\
\leq&~C_T \Bigg(\int_0^t\int_{\bbr^d} g^2(t-s,x-y) (1+\|\phi(s,y)\|^2_{L^2})\bone_{\{p=2\}}\,\dd(s,y)\\
& + \int_E j_N(z)\,\la(\dd z) \int_0^t\int_{\bbr^d} g^p(t-s,x-y)(1+\|\phi(s,y)\|_{L^p}^p) h(y)^{p-q}\,\dd(s,y)\Bigg)^{\frac{1}{p}}\\
\leq&~C_T\Bigg(\int_0^t\int_{\bbr^d} g^p(t-s,x-y)(1+\|\phi(s,y)\|_{L^p}^p) h(y)^{p-q}\,\dd(s,y)\Bigg)^{\frac{1}{p}}
\end{align*}
Next, by H\"older's inequality, we obtain
\begin{align*}
&~\left\|\int_0^t\int_{\bbr^d} G(t,x;s,y)\si(\phi(s,y))\,B_1^{N}(\dd s,\dd y)\right\|_{L^p}\\
\leq&~ C_T \left(\int_0^T\int_{\bbr^d} g(t,x)\,\dd(t,x)\right)^{1-\frac{1}{p}}\bbe\Bigg[\int_0^t\int_{\bbr^d} g(t-s,x-y)(1+|\phi(s,y)|^p)\\
&\times\left(|b(s,y)|+ \int_E |\ov\delta(s,y,z)|\bone_{\{|\delta(s,y,z)|\leq Nh(y)\}}\,\la(\dd z)\right)^p\bone_{\llb 0, \tau_N\rrb}(s)\,\dd(s,y)\Bigg]^{\frac{1}{p}}\\
\leq&~ C_T\left(\int_0^t\int_{\bbr^d} g(t-s,x-y)(1+\|\phi(s,y)\|_{L^p}^p)\left(\beta_N+\int_E j_N(z)\,\la(\dd z) (Nh(y))^{(1-q)_+}\right)^p\,\dd(s,y)\right)^{\frac{1}{p}}\\
\leq&~C_T\left(\int_0^t\int_{\bbr^d} g(t-s,x-y)(1+\|\phi(s,y)\|_{L^p}^p) h(y)^{p-q}\,\dd(s,y)\right)^{\frac{1}{p}},
\end{align*}
which completes the proof for the case $p\geq1$.

For $p<1$ we have by hypothesis $b_0, c\equiv0$. Therefore, 
\begin{align*} &~\left\|\int_0^t\int_{\bbr^d} G(t,x;s,y)\si(\phi(s,y))\,M^{N}(\dd s,\dd y)\right\|_{L^p}\\
\leq&~ C_T \bbe\left[\int_0^t\int_{\bbr^d}\int_E |g(t-s,x-y)\si(\phi(s,y))\delta(s,y,z)|^p\bone_{\{|\delta(s,y,z)|\leq Nh(y)\}}\bone_{\llb 0, \tau_N\rrb}(s)\,\pf(\dd s,\dd y,\dd z)\right]\\
\leq&~ C_T \int_0^t\int_{\bbr^d} g^p(t-s,x-y)(1+\|\phi(s,y)\|_{L^p}) h(y)^{p-q}\,\dd(s,y). \end{align*}
\halmos
\epr

We need two further preparatory results. The first one is classic. Henceforth, we denote by $\caln_d(\mu,\Sigma)$ the $d$-dimensional normal distribution with mean vector $\mu$ and covariance matrix $\Sigma$.
\blem\label{Gauss} If $X$ is an $\caln_1(0,\si^2)$-distributed random variable, then for every $p\in(-1,\infty)$ we have
\[ \bbe[|X|^p] = (2\si^2)^{p/2}\pi^{-1/2}\Ga\left(\textstyle\frac{1+p}{2}\right). \]
\elem

The second lemma determines the size of certain iterated integrals. It is proved by a straightforward induction argument, which is omitted.
\blem\label{itint}
Let $t\in\bbr_+$ and $a\in(-1,\infty)$. Then we have for every $n\in\bbn$ that ($t_0:=t$)
\[ \int_0^t\int_0^{t_1} \ldots\int_0^{t_{n-1}} \prod_{j=1}^n (t_{j-1}-t_j)^a \,\dd t_n\,\ldots\,\dd t_2\,\dd t_1 = \frac{\Ga(1+a)^n}{\Ga(1+(1+a)n)} t^{n(1+a)}.\]
\elem

We are now in the position to prove the existence of a solution to \eqref{SPDE-var} under the conditions of Theorem~\ref{main}.
\bpr[of Theorem~\ref{main}] (i)\quad We first prove that Equation~\eqref{SPDE-var} has a solution $Y^{(N)}$ in $B^p_\loc$ when $M$ is replaced by $M^{N}$ as defined in \eqref{MNK}. In order to do so, we choose the number $\eta$ in such a way that $\eta>d/q$ and $\eta(p-q)/2<1-d(p-1)/2$, which is possible by hypothesis A4. For reasons of readability, we do not index the subsequent processes with $N$ in this part of the proof, but only later when the dependence on $N$ matters. We define a Picard iteration scheme by $Y^0(t,x):=Y_0(t,x)$ and 
\beq Y^n(t,x):=Y_0(t,x)+\int_0^t\int_{\bbr^d} G(t,x;s,y)\si(Y^{n-1}(s,y))\,M^N(\dd s,\dd y),\quad (t,x)\in\bbr_+\times\bbr^d, \label{Yn}\eeq
for $n\in\bbn$. By Lemma~6.2 in \cite{Chong16} we can always choose a predictable version of $Y^n$. Then, since $Y_0\in B^p(\tau_N)$ and $\int_0^t\int_{\bbr^d} g^p(t-s,x-y)|y|^a\,\dd(s,y)<\infty$ for all $a\in\bbr_+$ and $p\in(0,1+2/d)$, Lemma~\ref{mom}(1) yields that $\|Y^n(t,x)\|_{L^p}<\infty$ for all $n\in\bbn$ and $(t,x)\in\bbr_+\times\bbr^d$. Now define $u^n:=Y^n-Y^{n-1}$ for $n\in\bbn$, which by Lemma~\ref{mom}(2) satisfies
\[ \|u^n(t,x)\|_{L^p} \leq C_T\left(\int_0^t \int_{\bbr^d} (g^p+g\bone_{\{p\geq1\}})(t-s,x-y)\|u^{n-1}(s,y)\|_{L^p}^{p^\ast} h(y)^{p-q}\,\dd(s,y)\right)^{\frac{1}{p^\ast}}\]
for all $n\geq2$. If we iterate this $n$ times, we obtain for all $(t,x)\in [0,T]\times[-R,R]^d$ (we abbreviate $g_p:=g^p+g\bone_{\{p\geq1\}}$)
\begin{align}
\|u^n(t,x)\|^{p^\ast}_{L^p} &\leq C^n_T \int_0^t \int_{\bbr^d} \ldots \int_0^{t_{n-1}}\int_{\bbr^d} g_p(t-t_1,x-x_1) \ldots g_p(t_{n-1}-t_n,x_{n-1}-x_n) \nonumber\\ 
&\quad\times (1+\|Y_0(s,y)\|_{L^p}^{p^\ast})h(x_1)^{p-q} \ldots h(x_n)^{p-q} \,\dd(t_n,x_n)\ldots\,\dd(t_1,x_1) \nonumber\\
&\leq C^n_T \int_0^t \ldots\int_0^{t_{n-1}} \int_{\bbr^d}\ldots\int_{\bbr^d} g_p(t-t_1,x_1) h(x-x_1)^{p-q} \ldots \nonumber\\
&\quad\times g_p(t_{n-1}-t_n, x_n) h(x-x_1-\ldots-x_n)^{p-q}\,\dd x_n\ldots\,\dd x_1\,\dd t_n\,\ldots\,\dd t_1. \label{help}
\end{align}
We take a closer look at the $n$ integrals with respect to $x_1, \ldots, x_n$. Define $\xi_i:=\sum_{j=1}^i x_j$ for $j=1,\ldots,n$ and $g_{p,1}:=g^p$ and $g_{p,2}:=g\bone_{\{p\geq1\}}$. Then, by H\"older's inequality, those integrals are bounded by
\begin{align} &~\sum_{l_1,\ldots,l_n = 1}^2 \int_{\bbr^d}\ldots\int_{\bbr^d} g_{p,l_1}(t-t_1,x_1) \ldots g_{p,l_n}(t_{n-1}-t_n, x_n) \prod_{i=1}^n h(x-\xi_i)^{p-q} \,\dd x_n\ldots\,\dd x_1\nonumber\\
\leq&~\sum_{l_1,\ldots,l_n = 1}^2 \prod_{i=1}^n \left(\int_{\bbr^d}\ldots\int_{\bbr^d} g_{p,l_1}(t-t_1,x_1)\ldots g_{p,l_n}(t_{n-1}-t_n,x_n) h(x-\xi_i)^{n(p-q)}\,\dd x_n\ldots\,\dd x_1\right)^{\frac{1}{n}}. \label{help2}\end{align}
We observe that $g(t,\cdot)$ is the density of the $\caln_d(0,2t\mathrm{I}_d)$-distribution, while $p^{d/2} (4\pi t)^{(p-1)d/2}g^p(t,\cdot)$ is that of the $\caln_d(0,2tp^{-1}\mathrm{I}_d)$-distribution. Here $\mathrm{I}_d$ is the $d$-dimensional identity matrix. Now let $X_1,\ldots,X_n$ be independent random variables such that $X_i$ has an $\caln_d(0,2(t_{i-1}-t_i)p^{-1}\mathrm{I}_d)$-distribution when $l_i=1$ and an $\caln_d(0,2(t_{i-1}-t_i)\mathrm{I}_d)$-distribution when $l_i=2$ (with $t_0:=t$). Then, with 
\[ \Xi_i:=\sum_{j=1}^i X_j,\quad \ga_i=\sum_{j\leq i, l_j=1} 2(t_{j-1}-t_j)p^{-1} + \sum_{j\leq i, l_j=2} 2(t_{j-1}-t_j),\quad i=1,\ldots,n,\] 
we use Lemma~\ref{Gauss} and the fact that $|x^{(k)}|\leq R$ and $\ga_i \leq 2(t-t_i)/(p\wedge1)\leq 2T/(p\wedge1)$ in order to derive
\begin{align*} &\int_{\bbr^d}\ldots\int_{\bbr^d} g_{p,l_1}(t-t_1,x_1)\ldots g_{p,l_n}(t_{n-1}-t_n,x_n) h(x-\xi_i)^{n(p-q)}\,\dd x_n\ldots\,\dd x_1\\
&\qquad=\left(\prod_{j\colon l_j=1} p^{-\frac{d}{2}}(4\pi (t_{j-1}-t_j))^{-\frac{d}{2}(p-1)}\right) \bbe[h(x-\Xi_i)^{n(p-q)}]\\
&\qquad \leq C_T^n \left(\prod_{j=1}^n (t_{j-1}-t_j)^{-\frac{d}{2}(p-1)}\right) \left(1+\bbe[|x-\Xi_i|^{n\eta(p-q)}]\right)\\
&\qquad\leq C_T^n \left(\prod_{j=1}^n (t_{j-1}-t_j)^{-\frac{d}{2}(p-1)}\right) \left(1+C_T^n \sup_{k=1,\ldots,d} \left(|x^{(k)}|^{n\eta(p-q)} + \bbe[|\Xi_i^{(k)}|^{n\eta(p-q)}]\right)\right)\\
&\qquad\leq C_T^n \left(\prod_{j=1}^n (t_{j-1}-t_j)^{-\frac{d}{2}(p-1)}\right) \left(1+C_T^n  \left(R^{n\eta(p-q)} + (2\ga_i)^{\frac{n\eta(p-q)}{2}}\pi^{-\frac{1}{2}}\Ga\left(\textstyle\frac{1+{n\eta(p-q)}}{2}\right) \right) \right)\\
&\qquad\leq C_T^n \left(\prod_{j=1}^n (t_{j-1}-t_j)^{-\frac{d}{2}(p-1)}\right)\Ga\left(\textstyle\frac{1+{n\eta(p-q)}}{2}\right).
\end{align*}
The last term no longer depends on $l_1,\ldots,l_n$ and $i$. Thus the last expression in \eqref{help2} is bounded by
\beq\label{help3} 2^n C^n_T \left(\prod_{j=1}^n (t_{j-1}-t_j)^{-\frac{d}{2}(p-1)}\right)\Ga\left(\textstyle\frac{1+{n\eta(p-q)}}{2}\right). \eeq
We insert this result back into \eqref{help} and apply Lemma~\ref{itint} to obtain
\begin{align} \|u^n(t,x)\|^{p^\ast}_{L^p} &\leq C^n_T \Ga\left(\textstyle\frac{1+{n\eta(p-q)}}{2}\right) \int_0^t \ldots\int_0^{t_{n-1}} \prod_{j=1}^n (t_{j-1}-t_j)^{-\frac{d}{2}(p-1)} \,\dd t_n\,\ldots\,\dd t_1,\nonumber\\
&\leq C^n_T  \Ga\left(\textstyle\frac{1+{n\eta(p-q)}}{2}\right) \frac{\Ga\left(1-\frac{d}{2}(p-1)\right)^n}{\Ga\left(1+(1-\frac{d}{2}(p-1))n\right)},\label{un} \end{align}
valid for all $(t,x)\in[0,T]\times[-R,R]^d$. Together with our choice of $\eta$ (see the beginning of the proof), it follows that 
\[ \sum_{n=1}^\infty \sup_{(t,x)\in[0,T]\times[-R,R]^d} \|u^n(t,x)\|_{L^p}<\infty,\]
implying that there exists some $Y\in B^p_\loc$ such that $\|Y^n-Y\|_{p,T,R}\to 0$ for every $T, R\in\bbr_+$ 
as $n\to\infty$. This $Y$ solves \eqref{SPDE-var} with $M^{N}$ instead of $M$. Indeed, the previous calculations actually show that for every $(t,x)\in\bbr_+\times\bbr^d$, $G(t,x;\cdot,\cdot)\si(Y^n)$ forms a Cauchy sequence with respect to the Daniell mean $\|\cdot\|_{M^{N},p}$. Therefore, it must converge to its limit $G(t,x;\cdot,\cdot)\si(Y)$ with respect to the Daniell mean, pointwise for $(t,x)$. In particular, the stochastic integrals with respect to $M^{N}$ converge to each other, showing that $Y$ is a solution. 

(ii)\quad As a next step we define for $n\in\bbn$ the operators 
\begin{align}\label{Jndef} J^{(1)}(\phi)&=J(\phi), &J^{(n)}(\phi)&:=J(J^{(n-1)}(\phi)),\\
	J^{(1)}_N(\phi)&=J_N(\phi), &J_N^{(n)}(\phi)&:=J_N(J_N^{(n-1)}(\phi)),\quad N\in\bbn,\quad n\geq2,\end{align}
hereby setting $J^{(n)}(\phi),J^{(n)}_N(\phi):=+\infty$ as soon as $J^{(n-1)}(\phi),J^{(n-1)}_N(\phi)=+\infty$. Then apart from the solution $Y=Y^{(N)}$ we found in (i) there exists no $\bar Y\in B^p_\loc$ which also satisfies $J_N(\bar Y)=\bar Y$ and for which there exists some $\bar Y_0\in B^p(\tau_N)$ such that $\bar Y^n:=J^{(n)}_N(\bar Y_0)$ converges to $\bar Y$ in $B^p_\loc(\tau_N)$ as $n\to\infty$. Indeed, by the same arguments as above one can show that $\|Y^n- \bar Y^n\|_{p,T,R}$ is bounded by the right-hand side of \eqref{un} to the power of $1/p^\ast$, possibly with another constant $C_T$. So taking the limit $n\to\infty$ proves $Y=\bar Y$.

(iii)\quad A solution to the original equation \eqref{SPDE-var} is now given by
\[ Y:= Y^{(1)}\bone_{\llb 0,\tau(1) \rrb}+\sum_{N=2}^\infty Y^{(N)} \bone_{\rrb \tau(N-1),\tau(N) \rrb},  \]
where $Y^{(N)}$ is the solution constructed in (i). By (ii), we have $Y^{(N)}\bone_{\llb 0,\tau(K)\rrb}=Y^{(K)}\bone_{\llb 0,\tau(K)\rrb}$ for all $N\in\bbn$ and $K=1,\ldots,N$.  Thus, as in the proof of Theorem~3.5 of \cite{Chong16}, one can verify that $Y$ indeed solves \eqref{SPDE-var}.
 \halmos\epr
 
 The uniqueness statement in part (ii) of the proof above can also be formulated for \eqref{SPDE-var}. 
 \bthm\label{uniqueness}
 The process $Y$ constructed in Theorem~\ref{main} is the unique solution to \eqref{SPDE-var} in the space of processes $\phi\in\tilde\calp$ for which there exist a sequence of $\bbf$-stopping times $(T_N)_{N\in\bbn}$ increasing to $+\infty$ a.s. and a process $\phi_0$ such that for arbitrary $T,R\in\bbr_+$ and $N\in\bbn$ we have $\phi_0\in B^p(T_N)$ and
 \[ \|(\phi-J^{(n)}(\phi_0))\bone_{\llb 0, T_N\rrb} \|_{p,T,R}  \to 0,\quad n\to\infty, \]
 where $J^{(n)}$ is defined via \eqref{Jndef}.
 \ethm
 
 Of course, this uniqueness result is quite weak: it does not say much about how other potential solutions to \eqref{SPDE-var} compare with the one we have constructed. Let us explain why we are not able to derive uniqueness in, say, $B^p_\loc(\tau(N))$ although the Picard iteration technique  --- via the Banach fixed point theorem --- usually yields existence and uniqueness at the same time. The reason is simply the following: on the one hand, the stopping times from Lemma~\ref{tauNK} allow us to obtain locally finite $L^p$-estimates. On the other hand, as one can see from Lemma~\ref{mom}, an extra factor $h(y)^{p-q}$ appears in these estimates. So for any $\phi\in\tilde\calp$ with $J_N(\phi)<\infty$ the right-hand side of \eqref{JN1} increases faster in $x$ than the input $\|\phi(t,x)\|_{L^p}$. Consequently, we cannot find a complete subspace of $B^p_\loc(\tau(N))$ containing $Y^0$ on which the operator $J_N$ is a self-map, which is a crucial assumption for the Banach fixed point theorem.
 
 Nevertheless, the next result demonstrates that the solution constructed in Theorem~\ref{main} is the ``natural'' one in terms of approximations. 
 \bthm\label{approx} Let $Y$ be the solution process to \eqref{SPDE-var} as constructed in Theorem~\ref{main} where the stopping times $\tau(N)$ are given by \eqref{tauN}. Furthermore, consider the following two ways of truncating the noise $M$:
 \benu
 \item $M_L(\dd t,\dd x):= b(t,x)\,\dd(t,x) + \rho(t,x)\,W(\dd t, \dd x) + \displaystyle\int_E \un\delta(t,x,z)\,(\pf-\qf)(\dd t, \dd x, \dd z)$ \\
 $~\qquad\qquad\qquad+\displaystyle\int_E \ov\delta(t,x,z)\bone_{\{|\delta(t,x,z)|\leq L\}}\,\pf(\dd t,\dd x,\dd z),\quad L\in \bbn$,
 \item $M_L(\dd t, \dd x) := b(t,x)\,\dd(t,x) + \rho(t,x)\,W(\dd t, \dd x) + \displaystyle\int_E \un\delta(t,x,z)\,(\pf-\qf)(\dd t, \dd x, \dd z)$ \\
 $~\qquad\qquad\qquad+\displaystyle\int_E \ov\delta(t,x,z)\bone_{[-L,L]^d}(x)\,\pf(\dd t,\dd x,\dd z),\quad L\in \bbn$.
 \eenu
 In both cases, if $Y_L$ denotes the unique solution to \eqref{SPDE-var} with $M$ replaced by $M_L$ that belongs to $B^p(\tau(N))$ for all $N\in\bbn$ (see Theorems~3.1 and 3.5 in \cite{Chong16}), then we have for all $N\in\bbn$ and $T,R\in\bbr_+$ 
 \[ \lim_{L\to\infty}\|(Y_L-Y)\bone_{\llb 0, \tau(N)\rrb}\|_{p,T,R}=0.\]
 \ethm
 
 From another point of view, Theorem~\ref{approx} paves the way for simulating from the solution of \eqref{SPDE-var}. In fact, in \cite{Chen16} different methods are suggested for the simulation of \eqref{SPDE-var} with the noise $M_L$ as given in the first part of the previous theorem. The approximations in that paper were shown to be convergent in an $L^p$-sense. Thus, letting $L$ increase in parallel, these approximations will then converge to the solution of \eqref{SPDE-var} with the untruncated noise $M$, at least in $B^p_\loc(\tau(N))$.
 
\bpr[of Theorem~\ref{approx}] (i)~We first prove the case (1). Let $N$, $T$ and $R$ be fixed. Then the arguments in the proof of Theorem~\ref{main} reveal that
\[ \lim_{n\to\infty} \|(Y-Y^n)\bone_{\llb 0,\tau(N)\rrb}\|_{p,T,R} = \lim_{n\to\infty} \sup_{L\in\bbn} \|(Y_L-Y^n_L)\bone_{\llb 0,\tau(N)\rrb}\|_{p,T,R}=0,  \]
where $Y^n$ is given by \eqref{Yn} and $Y^n_L$ is defined in the same way but with $M^N_L$ instead of $M^N$, and $M^N_L$ is the measure specified through \eqref{MNK} with $Nh(x)$ replaced by $Nh(x)\wedge L$. Therefore, it suffices to prove that for every $n\in\bbn$ we have $\|(Y^n-Y^n_L)\|_{p,T,R}\to0$ as $L\to\infty$. To this end, we observe that 
\begin{align} Y^n(t,x)-Y^n_L(t,x) &= \int_0^t\int_{\bbr^d} G(t,x;s,y)\big(\si(Y^{n-1}(s,y)) -\si(Y^{n-1}_L(s,y))\big)\,M^N_L(\dd s,\dd y) \nonumber\\
&\quad+ \int_0^t\int_{\bbr^d} G(t,x;s,y)\si(Y^{n-1}(s,y))\,(M^N-M^N_L)(\dd s,\dd y)\nonumber\\
&=:v^n(t,x)+w^n(t,x). \label{eq1} \end{align}
Furthermore, since $n$ is fixed and the function $|x|\mapsto\int_0^T\int_{\bbr^d} |g(s,x-y)|^p|y|^a\,\dd(s,y)$ grows at most polynomially in $|x|$ for any $a\in\bbr^+$, it follows that $\|Y^{n-1}(t,x)\|^p_{L^p}\leq C_T (1+|x|)^{m/2}$ for all $(t,x)\in[0,T]\times\bbr^d$ when $m\in\bbr_+$ is chosen large enough (we may assume that $m/2\geq\eta(p-q)$). Therefore, recalling the notation $g_p=g^p+g\bone_{\{p\geq1\}}$, we obtain for $p\geq1$
\begin{align*} \|w^n(t,x)\|^p_{L^p} 
&\leq  \Bigg\|\int_0^t\int_{\bbr^d}\int_E G(t,x;s,y)\si(Y^{n-1}(s,y))\delta(s,y,z)\bone_{\{L< |\delta(s,y,z)|\leq Nh(y) \}}\\
&\quad\times\bone_{\llb 0,\tau(N)\rrb}(s)\,(\pf-\qf)(\dd s,\dd y,\dd z)\Bigg\|^p_{L^p} \\
&\quad+\Bigg\| \int_0^t\int_{\bbr^d}\int_E G(t,x;s,y)\si(Y^{n-1}(s,y))\delta(s,y,z)\bone_{\{L< |\delta(s,y,z)|\leq Nh(y) \}}\\
&\quad\times\bone_{\llb 0,\tau(N)\rrb}(s)\,\qf(\dd s,\dd y,\dd z)\Bigg\|^p_{L^p}\\
&\leq C_T\Bigg(\int_0^t \int_{\bbr^d}g^p(t-s,x-y) \|\si(Y^{n-1}(s,y))\|^p_{L^p}  (Nh(y))^{p-q}\bone_{\{Nh(y)>L\}}\,\dd(s,y)\\
&\quad+ \int_0^t\int_{\bbr^d} g(t-s,x-y)\|\si(Y^{n-1}(s,y))\|^p_{L^p}  (Nh(y))^{p-q}\bone_{\{Nh(y)>L\}}\,\dd(s,y)\\
&\quad\times \left(\int_0^t\int_{\bbr^d} g^p(s,y)\,\dd(s,y)\right)^{p-1}\Bigg)\\
&\leq C_T \int_0^t\int_{\bbr^d} g_p(s,x-y) (1+|y|)^m \bone_{\big\{|y|>\left(\frac{L}{N}-1\right)^{1/\eta}\big\}}\,\dd(s,y).
\end{align*}
Using slightly modified calculations, one can derive the final bound also for $p<1$ and show that
\[ \|v^n(t,x)\|_{L^p}^{p^\ast} \leq C_T \int_0^t\int_{\bbr^d} g_p(t-s,x-y)\|Y^{n-1}(s,y)-Y^{n-1}_L(s,y)\|_{L^p}^{p^\ast}h(y)^{p-q}\,\dd(s,y) \]
with a constant $C_T$ independent of $L$.
Plugging these estimates into \eqref{eq1}, we obtain by iteration for all $(t,x)\in[0,T]\times[-R,R]^d$ ($\theta_i:=\sum_{j=1}^i t_i$ and $\xi_i:=\sum_{j=1}^i x_i$)
\begin{align*}
&\|Y^n(t,x)-Y^n_L(t,x)\|_{L^p}^{p^\ast}\\
&\qquad \leq C_T\int_0^t\int_{\bbr^d} g_p(t-s,x-y) \|Y^{n-1}(s,y)-Y^{n-1}_L(s,y)\|_{L^p}^{p^\ast}(1+|y|)^m\,\dd(s,y)\\
&\qquad\quad +C_T\int_0^t\int_{\bbr^d} g_p(t-s,x-y) (1+|y|)^m \bone_{\big\{|y|>\left(\frac{L}{N}-1\right)^{1/\eta}\big\}}\,\dd(s,y)\\
&\qquad\leq\sum_{j=1}^n C_T^j  \int_0^t\int_{\bbr^d} \ldots \int_0^{t_{j-1}}\int_{\bbr^d} g_p(t-t_1,x-x_1)\ldots g_p(t_{j-1}-t_j,x_{j-1}-x_j)\\
&\qquad\quad\times(1+|x_1|)^m\ldots (1+|x_j|)^m \bone_{\big\{|x_j|>\left(\frac{L}{N}-1\right)^{1/\eta}\big\}}\,\dd(t_j,x_j)\ldots\,\dd(t_1,x_1)\\
&\qquad =\sum_{j=1}^n C_T^j  \int_0^t\int_{\bbr^d} \ldots \int_0^{t-\theta_{j-1}} \int_{\bbr^d} g_p(t_1,x_1)\ldots g_p(t_j,x_j)\\
&\qquad\quad\times(1+|x-\xi_1|)^m \ldots (1+|x-\xi_j|)^m \bone_{\big\{|x-\xi_j|>\left(\frac{L}{N}-1\right)^{1/\eta}\big\}}\,\dd(t_j,x_j)\ldots\,\dd(t_1,x_1)\\
&\qquad\leq\sum_{j=1}^n C_T^j  \int_0^T\int_{\bbr^d} \ldots \int_0^T \int_{\bbr^d} g_p(t_1,x_1)\ldots g_p(t_j,x_j)(1+d^{1/2}R+|\xi_1|)^m\\
&\qquad\quad\times \ldots \times(1+d^{1/2}R+|\xi_j|)^m \bone_{\big\{|\xi_j|>\left(\frac{L}{N}-1\right)^{1/\eta}-d^{1/2}R\big\}}\,\dd(t_j,x_j)\ldots\,\dd(t_1,x_1).
\end{align*}
The integrals in the last line do not depend on $(t,x)$ and would be defined even without the indicator function. Hence they converge to $0$ as $L\to\infty$ by dominated convergence.

(ii)~If $M_L$ is defined as in (2), we first notice that for every $L,N\in\bbn$ there exists $N^\prime\in\bbn$ such that $\tau(N)\leq \tau(N^\prime,L)$ a.s. where
\[ \tau(N,L):=\inf\left\{T\in\bbr_+\colon \int_0^T\int_{[-L,L]^d}\int_E \bone_{\{|\bar\delta(t,x,z)>N|\}}\,\pf(\dd t,\dd x,\dd z) \neq 0\right\},\quad N,L\in\bbn. \]
Therefore, Theorem~3.5 in \cite{Chong16} which states that there exists a unique solution to \eqref{SPDE-var} with noise $M_L$ that belongs to $B^p(\tau(N,L))$ for all $N\in\bbn$ automatically implies that this solution is also the unique one that belongs to $B^p(\tau(N))$ for all $N\in\bbn$. The actual claim is proved in basically the same way as in (1), except that in the moment estimate of $w^n(t,x)$ one has to replace the indicator $\bone_{\{Nh(y)>L\}}$ by $\bone_{\bbr^d\setminus[-L,L]^d}(y)$, which obviously does not effect the final convergence result.
\halmos
\epr

Until now, the solution $Y$ to \eqref{SPDE-var} that we have constructed in Theorem~\ref{main} only has finite $p$-th moments, locally uniformly in space and locally uniformly in time until $T\wedge \tau(N)$ for any $T$ and $N$. But does $Y$ possess any finite moments until $T\wedge\tau_N$ where $\tau_N$ is the original stopping time from hypothesis M1? If $M$ is a L\'evy basis, this would be the question whether $Y$ has finite moments up to any fixed time $T$ without stopping. Furthermore, under which conditions do they remain bounded and not blow up in space? The next theorem provides a sufficient condition for these statements to be true.
\bthm\label{moment}
If additionally to the assumptions of Theorem~\ref{main} we have that 
\[ |\si(x)|\leq C(1+|x|^\ga),\quad x\in\bbr,\]
for some $C\in\bbr_+$ and $\ga\in [0,q/p]$, then the solution $Y$ to \eqref{SPDE-var} as constructed in Theorem~\ref{main} belongs to $B^q(\tau_N)$ for all $N\in\bbn$.
\ethm

We need an analogue of Lemma~\ref{mom} for moments of order $q$.
\blem\label{mom2} Let $J$ be the integral operator defined in \eqref{Jdef}. For all $K\in\bbn$ and $T\in\bbr_+$ we can find a constant $C_T\in\bbr_+$ 
such that for all $\phi\in\tilde\calp$ and $(t,x)\in[0,T]\times\bbr^d$ we have
\begin{align*} &\|J(\phi)(t,x)\bone_{\llb 0, \tau_K\rrb}(t)\|_{L^q} \leq \|Y_0(t,x)\bone_{\llb 0, \tau_K\rrb}(t)\|_{L^q} \\
	&\quad\quad+ C_T \left(\int_0^t\int_{\bbr^d} g^p(t-s,x-y)\|\si(\phi(s,y))\bone_{\llb 0, \tau_K\rrb}(s)\|^{p^\ast}_{L^p}\,\dd(s,y)\right)^{q_\ast/p} \\
	&\quad\quad+C_T \left(\int_0^t\int_{\bbr^d} g^q(t-s,x-y)(1+\|\phi(s,y)\bone_{\llb 0, \tau_K\rrb}(s)\|^{q^\ast}_{L^q})\,\dd(s,y)\right)^{1/q^\ast}\\
	&\quad\quad+ C_T \left(\int_0^t\int_{\bbr^d} g(t-s,x-y)\bone_{\{p\geq1,q\geq 1\}}(1+\|\phi(s,y)\bone_{\llb 0, \tau_K\rrb}(s)\|^{q}_{L^q})\,\dd(s,y)\right)^{1/q}\\
	&\quad\quad+ C_T \left(\int_0^t\int_{\bbr^d} g(t-s,x-y)\bone_{\{p\geq1,q< 1\}}\|\si(\phi(s,y))\bone_{\llb 0, \tau_K\rrb}(s)\|^p_{L^p}\,\dd(s,y)\right)^{q/p}. \end{align*}
Moreover, if $J_N$ is the integral operator defined in Lemma~\ref{mom}, then the previous estimates also hold for $J_N(\phi)$ and the constant $C_T$ does not depend on $N$.
\elem
\bpr In principle the proof follows the same line of reasoning as Lemma~\ref{mom}. If $p,q\geq1$, the idea is to split the It\^o basis $M$ into three parts:
\begin{align*} M(\dd t, \dd x) &=\left(b(t,x)+\int_E \ov\delta(t,x,z)\,\la(\dd z)\right)\,\dd (t,x) +\int_E \ov\delta(t,x,z)\,(\pf-\qf)(\dd t,\dd x,\dd z)\,\dd (t,x) \\
	&\quad+\left(\rho(t,x)\,W(\dd t,\dd x) + \int_E \un\delta(t,x,z)\,(\pf-\qf)(\dd t,\dd x,\dd z)\right). \end{align*}
The $q$-th moments of the integrals against the first two parts can be estimated as in Lemma~\ref{mom}. Since $q$ replaces $p$, the factor $h(y)^{p-q}$ can be omitted throughout. For the integral against the third part in the decomposition above, its $L^q$-norm is bounded by its $L^p$-norm, which can be treated as in Lemma~\ref{mom}. Again, the factor $h(y)^{p-q}$ is not needed because the jump sizes are at most $1$. Next, if $p\geq1$ but $q<1$, we can proceed in the same way except that we replace the first and second part in the decomposition of $M$ by
\[ b(t,x)\,\dd(t,x)\quad\text{and}\quad \int_E \ov\delta(t,x,z)\,\pf(\dd t,\dd x,\dd z), \]
respectively. Finally, if $p,q<1$, we have that 
\[ M(\dd t,\dd x) = \int_E \un\delta(t,x,z)\,\pf(\dd t,\dd x,\dd z) + \int_E \ov\delta(t,x,z)\,\pf(\dd t,\dd x,\dd z). \] 
The result then follows in a similar way if we switch to the $p$-th moment for the first term. The statement regarding $J_N$ is evident because the estimates above hold for all measures with only a subset of jumps of $M$.\halmos
\epr

\bpr[of Theorem~\ref{moment}]  Let $N\in\bbn$ and $Y^{n,N}$ be the process defined in \eqref{Yn}. Then, because we have $|\si(x)|\leq C(1+|x|^\ga)$ with $\ga\leq q/p$, Lemma~\ref{mom2} implies that for every $K\in\bbn$ we have
\begin{align*}
\|Y^{n,N}(t,x)\bone_{\llb 0,\tau_K\rrb}(t)\|_{L^q}^{q^\ast} &\leq C_T\Bigg(1+\int_0^t\int_{\bbr^d} g^p(t-s,x-y) \|Y^{n-1,N}(s,y)\bone_{\llb 0,\tau_K\rrb}(s)\|_{L^q}^{q^\ast}\,\dd(s,y) \\
&\quad+\int_0^t\int_{\bbr^d} g^q(t-s,x-y) \|Y^{n-1,N}(s,y)\bone_{\llb 0,\tau_K\rrb}(s)\|^{q^\ast}_{L^q}\,\dd(s,y)\\
&\quad+\int_0^t\int_{\bbr^d} g(t-s,x-y) \bone_{\{p\geq1\}}\|Y^{n-1,N}(s,y)\bone_{\llb 0,\tau_K\rrb}(s)\|^{q^\ast}_{L^q}\,\dd(s,y)\Bigg)
\end{align*}
with a constant $C_T$ independent of $N$ and $n$. Since $\|Y^{0,N}\bone_{\llb 0,\tau_K\rrb}\|_{q,T}=\|Y_0\bone_{\llb 0,\tau_K\rrb}\|_{q,T}<\infty$, we can apply Lemma~6.4 in \cite{Chong16} to deduce 
\[ \sup_{n,N\in\bbn} \|Y^{n,N}\bone_{\llb 0,\tau_K\rrb}\|_{q,T} <\infty. \]
Since $Y^{(N)}$ is the limit in $B^p_\loc$ of $Y^{n,N}$ as $n\to\infty$ and $Y$ is piecewise equal to $Y^{(N)}$, the assertion follows and $Y$ belongs to $B^q(\tau_K)$ for all $K\in\bbn$.
\halmos
\epr

How do the results and methods described above extend to more general equations than \eqref{SPDE-var}, for example, when $G$ satisfies A3 with some other kernel $g$ than the heat kernel? Immediately one notices the following difference between the case $q=p$ and $q<p$: while in the former case one obtains existence and uniqueness as soon as
\beq\label{intass} \int_0^T\int_{\bbr^d} |g(t,x)|^p\,\dd(s,y)<\infty \eeq
for all $T>0$ (see \cite{Chong16}, Theorem~3.1), we need to put a much stronger assumption in the latter case. Namely, since an extra factor $h(y)^{p-q}$ appears in each iteration step in the proof of Theorem~\ref{main}, the kernel $g$ must decay faster in $|x|$ than any polynomial. Thus, our methods will \emph{not} work merely under an integrability assumption like \eqref{intass}. For example, fractional equations as considered, for example, in \cite{Balan14}, \cite{Peszat07} or \cite{Wu12} are outside the scope of this paper when $q<p$. However, our results can be extended to quite general parabolic stochastic PDEs. More precisely, consider 
the partial differential operator given by
\beq\label{PDEoperator} \calL = \partial_t - \sum_{|\al|\leq 2m} a_\al(t,x)\partial^\alpha,\quad (t,x)\in\bbr_+\times\bbr^d, \eeq
with $m\in\bbn$ and suitable bounded continuous functions $a_\al$. Under regularity conditions (see \cite{Peszat07}, Theorem~2.6, and \cite{Eidelman98}, Theorem~VI.2), the operator \eqref{PDEoperator} admits a fundamental solution $G(t,x;s,y)$ satisfying
\beq
|G(t,x;s,y)| \leq C_T g_\calL(t-s,x-y),\quad g_\calL(t,x):=\frac{1}{t^{d/(2m)}}\exp\left( -\La\frac{|x|^{(2m)/(2m-1)}}{t^{1/(2m-1)}} \right)\bone_{(0,\infty)}(t),\label{genheat}\\
\eeq
for $(t,x),(s,y)\in\bbr_+\times\bbr^d$ with some strictly positive constants $\La$ and $C_T$ independent of $m$. With this at hand, we can now extend the previous results to incorporate the case where $G$ satisfies \eqref{genheat}. In fact, we can go one step further and consider kernels $G$ bounded by generalized Gaussian densities \citep{Gomez98} of the form
\beq\label{genGaus} g_{\rho,\tau,\La}(t,x) := K(\rho,\tau,\La) t^{-\tau d/\rho}\ee^{-\La |x|^\rho/t^\tau} \bone_{(0,\infty)}(t),\quad (t,x)\in\bbr_+\times\bbr^d, \eeq
for some parameters $\rho,\tau,\La>0$ and a constant $K(\rho,\tau,\La)$ chosen in such a way that $g_{\rho,\tau,\La}(t,\cdot)$ becomes a probability density function on $\bbr^d$ for every fixed $t$. 

\bthm\label{extension}
Assume the same hypotheses as in Theorem~\ref{main} except that in A3 the heat kernel $g$ is replaced by the function $g_{\rho,\tau,\La}$ from \eqref{genGaus} and A4 by the requirement 
\beq\label{A4new} 0<p<1+\frac{\rho}{\tau d}\quad\text{and}\quad \frac{p}{1+\tau(1+\frac{\rho}{\tau d}-p)} < q \leq p. \eeq
Then the assertions of Theorems~\ref{main}, \ref{uniqueness}, \ref{approx} and \ref{moment} continue to hold. 

In particular, in the case of \eqref{genheat}, we have $\rho=(2m)/(2m-1)$ and $\tau=1/(2m-1)$, so \eqref{A4new} becomes
\beq\label{A4new2} 0<p<1+\frac{2m}{d}\quad\text{and}\quad \frac{p}{1+(1+\frac{2m}{d}-p)/(2m-1)} < q \leq p. \eeq
\ethm

\bpr
The only step that has to be modified in the proofs when the heat kernel $g$ is replaced by $g_{\rho,\tau,\La}$ is the estimation of the integrals in \eqref{help2}. Because
\[ g_{\rho,\tau,\La}^p(t,x) = \frac{K(\rho,\tau,\La)^p}{K(\rho,\tau,p\La)} t^{-(p-1)\tau d/\rho}g_{\rho,\tau,\La p}(t,x),\quad(t,x)\in\bbr_+\times\bbr^d, \]
and
\begin{align*} \int_{\bbr^d}  g_{\rho,\tau,\La p}(t,x) |x|^{n\eta(p-q)} \,\dd x& = \frac{K(\rho,\tau,p\La)}{t^{\tau d/\rho}} \int_0^\infty \ee^{-p\La r^\rho/t^\tau} r^{n\eta(p-q)}r^{d-1}\,\dd r\\ &= \frac{K(\rho,\tau,p\La)}{\rho} t^{n\eta(p-q)/\rho} (p\La)^{-(d+n\eta(p-q))/\rho} \Gamma\left(\textstyle\frac{d+n\eta(p-q)}{\rho}\right), 
	\end{align*}
we obtain in the end instead of \eqref{help3} the bound
\[ C_T^n\left(\prod_{j=1}^n (t_{j-1}-t_j)^{-(p-1)\tau d/\rho} \right)\Gamma\left(\textstyle\frac{d+n\eta(p-q)}{\rho}\right), \]
and instead of \eqref{un} the estimate
\[ C^n_T\Gamma\left(\textstyle\frac{d+n\eta(p-q)}{\rho}\right) \frac{\Ga\left(1-(p-1)\tau d/\rho\right)^n}{\Ga\left(1+(1-(p-1)\tau d/\rho)n\right)}. \]
By \eqref{A4new} we can choose $\eta>d/q$ in such a way that $\eta(p-q)/\rho<1-(p-1)\tau d/\rho$ such that the previous bound converges to $0$ as $n$ tends to infinity.\halmos
\epr

\subsection*{Acknowledgement}
I wish to thank Jean Jacod for his constructive comments and suggestions. 

\addcontentsline{toc}{section}{References}
\bibliographystyle{plainnat}
\bibliography{bib-SPDE-stable}
\end{document}